\documentclass[a4paper,twoside,11pt]{article}

\usepackage{a4wide, fancyhdr, amsmath, amssymb, mathtools, yfonts}
\usepackage{mathrsfs}
\usepackage{graphicx}
\usepackage{tikz}
\usepackage[all]{xy}
\usepackage[utf8]{inputenc}
\usepackage{amsthm}
\usepackage[english]{babel}
\usepackage{chngcntr}
\usepackage{ifthen}
\usepackage{calc}
\usepackage{hyperref}
\usepackage{authblk}
\numberwithin{equation}{section}
\newcommand{\Q}{\mathbb{Q}}

\newtheorem{lemma}{Lemma}[section]
\newtheorem{theorem}[lemma]{Theorem}
\newtheorem{proposition}[lemma]{Proposition}
\newtheorem{corollary}[lemma]{Corollary}
\newtheorem{definition}[lemma]{Definition}

\title{\vspace{ - \baselineskip} \sffamily\bfseries The size of arboreal images, I: exponential lower bounds for PCF and unicritical polynomials}
\author[1,2]{Carlo Pagano\thanks{Vivatsgasse 7, 53111 Bonn, Germany, carlein90@gmail.com}}
\affil[1]{Max Planck Institute for Mathematics, Bonn}
\affil[2]{University of Glasgow, Glasgow}

\date{\today}

\begin{document}
\maketitle

\begin{abstract}
Let $f$ be a polynomial over a global field $K$. For each $\alpha$ in $K$ and $N$ in $\mathbb{Z}_{\geq 0}$ denote by $K_N(f,\alpha)$ the \emph{arboreal field} $K(f^{-N}(\alpha))$ and by $D_N(f,\alpha)$ its degree over $K$.

It is conjectured that $D_N(f,\alpha)$ should grow as a double exponential function of $N$, unless $f$ is \emph{post-critically finite} (PCF), in which case there are examples like $D_N(x^2,\alpha) \leq 4^{N}$. There is evidence conditionally on Vojta's conjecture. However, before the present work, no unconditional non-trivial lower bound was known for \emph{post-critically infinite} $f$. In the case $f$ is PCF, no non-trivial lower bound was known, not even under Vojta's conjecture. 

In this paper we give two simple methods that turn the finiteness of the critical orbit into an exploitable feature, also in the post-critically infinite case. First, assuming GRH for number fields, we establish for all PCF polynomials $f$ of degree at least $2$ and all $\alpha$ outside of the critical orbits of $f$, the existence of a positive constant $c(f,\alpha)$ such that
$$D_N(f,\alpha) \geq \text{exp}(c(f,\alpha)N),
$$
for all $N$, which is sharp up to, possibly, improve the constant $c(f,\alpha)$. 

Second, we show unconditionally that if $f$ is post-critically infinite over any number field $K$ and unicritical, then, for each $\alpha$ in $K$, there exists a positive constant $c(f,\alpha)$ such that
$$D_N(f,\alpha) \geq \text{exp}(c(f,\alpha)N),
$$
for all $N$. The main input here is to work modulo a suitably chosen prime and use a construction available for PCF unicritical polynomials with periodic critical orbit. 
%of the field generated by the
%first backward iterations $\beta \mapsto\{\gamma:~ f(\gamma) = \beta\}$, which start from some element $\alpha \in K$.
\end{abstract}

\maketitle

\section{Introduction}
Arboreal Galois groups are among the central objects of study in modern arithmetic dynamics. We refer to \cite[Section 5]{BIJMST} for a comprehensive list of results and open problems on the subject. Arboreal Galois groups are typically obtained by looking at the graph of pre-images 
$$T_{\infty}(f,\alpha):=\{f^{-N}(\alpha)\}_{N \in \mathbb{Z}_{\geq 0}}
$$
of an element $\alpha$ of a field $K$ under a separable polynomial $f \in K[x]$: the edges are defined by connecting $\beta$ to $f(\beta)$.  Under the extra assumption that $\alpha$ is not a periodic point and is not in the orbit of a critical point of $f$, one sees that $T_{\infty}(f,\alpha)$ is naturally an infinite degree $d$ tree rooted at $\alpha$. \footnote{Unless explicitly stated we won't make this extra assumption in what follows. Also we notice that it is enough to assume that $\alpha$ is not of the form $f^n(\beta)$ for $\beta$ a critical point and $n$ a \emph{positive} integer.} The absolute Galois group $G_K$ acts by graph-preserving symmetries on $T_{\infty}(f,\alpha)$ and gives rise to the \emph{arboreal Galois group}
$$G_{\infty}(f,\alpha):=\text{Gal}(K_{\infty}(f,\alpha)/K),
$$
where $K_{\infty}(f,\alpha)$ is the extension of $K$ generated by $T_{\infty}(f,\alpha)$, and more generally $K_N(f,\alpha):=K(f^{-N}(\alpha))$ for each $N$ in $\mathbb{Z}_{\geq 0}$, is the splitting field over $K$ of the polynomial $f^N-\alpha$. In this paper $f^N$ will always mean the $N$-th iterate of $f$ under composition.  

Let us now consider the case that $K$ is a global field. It is generally expected that arboreal Galois groups over $K$ should be \emph{large}. The first results in this direction are due to Odoni \cite{Odoni1, Odoni2}. Subsequently Cremona \cite{Cremona} considered the polynomial $x^2+1$ in detail and conjectured that for $K=\Q$ the arboreal Galois groups $G_N(x^2+1,0)$ are as large as possible. This was proved by Stoll \cite{Stoll}. Stoll's ingenious argument has been recently extended in \cite{FPC}, which considered more general open subgroups of the entire automorphism group of the tree $T_{\infty}(f,\alpha)$, which we denote by $\Omega_{\infty}(f,\alpha)$. Examples of surjectivity for higher degree polynomials can be found in \cite{Benedetto: surjectivity, Looper, Kadet, Spechter}. Some further explicit families of unicritical polynomials yielding an open image have been exihibited in   \cite{Bush--Hindes--Looper}, using the Chabauty--Coleman method and the Mordell--Weil sieve. 

A precise conjecture, reflecting the general expectation on arboreal images over global fields, is due to Jones, who conjectured \cite[Conjecture 3.11]{Jones} that for quadratic polynomials $f$ over a number field, with $\alpha$ not $f$-periodic and outside of the strict critical orbit of $f$, the group $G_{\infty}(f,\alpha)$ is an open subgroup of $\Omega_{\infty}(f,\alpha)$ precisely when $f$ is \emph{post-critically infinite}. We recall the definition of post-critically (in)finite. 

\begin{definition} \label{PCF}
Let $K$ be a field and $d$ be in $\mathbb{Z}_{\geq 2}$. Let $f \in K[x]$ be a monic polynomial of degree $d$. We say that $f$ is \emph{post-critically finite} (abbreviated \emph{PCF}) in case for each $\gamma \in \bar{K}$ with $f'(\gamma)=0$, we have that
$$\#\{f^n(\gamma) \}_{n \in \mathbb{Z}_{\geq 0}}<\infty. 
$$
If $f$ is not post-critically finite, we say that it is \emph{post-critically infinite}.

\end{definition}
We call each of the orbits $\{f^n(\gamma) \}_{n \in \mathbb{Z}_{\geq 0}}$ (respectively $\{f^n(\gamma) \}_{n \in \mathbb{Z}_{\geq 1}}$), with $\gamma$ a zero of $f'(x)$, a critical orbit of $f$ (respectively a strict critical orbit). The notion of PCF can be extended to all rational maps and in a precise sense PCF rational maps over a number field are not abundant: \cite{attracting cycles} shows that, aside from the set of Latt\`es maps, they form a set of bounded height in the moduli space of rational maps. 

Jones' conjecture would in particular imply that the finite arboreal Galois groups attached to post-critically infinite quadratic polynomials are very large in size. Namely that, as long as $(f,\alpha)$ obeys the above assumptions, for
$$G_N(f,\alpha):=\text{Gal}(K_N(f,\alpha)/K)
$$
one has that there exists a positive constant $c(f,\alpha)$ such that
$$\#G_N(f,\alpha) \geq c(f,\alpha) \cdot 2^{2^N}
$$
as $N$ varies among the positive integers. \footnote{Dropping the assumptions on $\alpha$, but keeping that $f$ is post-critically infinite, would still give a double exponential lower bound, with an adaptation very similar to the one that we will give at the end of the proof of Theorem \ref{exponential:GRH}.}Under Vojta's conjectures for number fields there is substantial evidence for Jones' conjecture \cite{Hindes}, and over function fields one can, in some cases, use these techniques to prove unconditionally that the representation is surjective \cite{FM}. Under the abc conjecture \cite[Theorem 3]{Win} shows Jones' conjecture in the case of eventually stable pairs $(f,\alpha)$. Here eventually stable means that $f^N-\alpha$ factorizes in a number of irreducibles that is eventually constant in $N$. Similar results were obtained in \cite{Gratton--Nguyen--Tucker}. Some unconditional evidence for Jones' conjecture over number fields has been provided in the recent work \cite{FP}.  

For the PCF case, proving lower bounds for $\#G_N(f,\alpha)$ seems at first a more elusive problem, already in the quadratic case and even if one assumes strong diophantine conjectures. Indeed, when $f$ is PCF, all the standard sub-extensions that one typically constructs out of the critical orbit of $f$, amount to a finite multi-quadratic subextension of $K_{\infty}(f,\alpha)/K$.  

In the present paper we give two simple methods that turn the finiteness of the critical orbit into a \emph{feature} that allows to prove lower bounds for $\#G_N(f,\alpha)$, both in the PCF and, perhaps surprisingly, in the post-critically infinite case. In the PCF case, our lower bound for $\#G_N(f,\alpha)$ are basically sharp. For the first method, however, instead of relying on diophantine conjectures we have to rely on conjectures from analytic number theory. In contrast our second method leads to an unconditional result. We recall that the pair $(f,\alpha)$ is said to be \emph{exceptional} in case $T_{\infty}(f,\alpha)$ is a finite set. A sufficient condition to be non exceptional is that $\alpha$ is not in any critical orbit of $f$, which excludes only finitely many values in case $f$ is PCF. Our first main theorem is as follows.

\begin{theorem} \label{exponential:GRH}
Let $K$ be a number field and let $f$ be a post-critically finite polynomial over $K$ of degree at least $2$. Let $\alpha$ be in $K$ and suppose $(f,\alpha)$ is not exceptional. Assume that GRH holds. Then there exists a positive constant $c(f,\alpha)$ such that
$$\#G_N(f,\alpha) \geq \textup{exp}(c(f,\alpha) \cdot N), 
$$
for all $N$ in $\mathbb{Z}_{\geq 1}$. 
 \end{theorem}
We recall that one has examples such as
$$\#G_N(x^2,\alpha) \leq 4^{N} \ \text{and} \ \#G_N(x^2-2,0) \leq 2^N,
$$
hence showing that our exponential lower bound for PCF is sharp, apart from possibly improving on the constant $c(f,\alpha)$. That said, exponential behavior might likely be very rare even among PCF and it might be that for arboreal Galois groups that are sufficiently non-abelian a much better lower bound holds. However, currently we do not know how to achieve this. Abelian arboreal Galois groups are the subject of a very recent conjecture \cite{Andrews--Petsche}, where the conjecture has been also proved for quadratic stable polynomials over $\Q$ \cite[Theorem 1]{Andrews--Petsche}. Later on, the conjecture has been settled for general quadratic polynomials over $\Q$ in \cite[Theorem 3.12]{FP}. In this work it was also shown that for a quadratic polynomial $f$ over a number field $K$ the associated arboreal representation is abelian for some $\alpha$ in $K$ only if $f$ is PCF \cite[Theorem 3.8]{FP}. 

In the proof of Theorem \ref{exponential:GRH}, we actually give as a lower bound $c(\epsilon, f,\alpha) \cdot (d^{\frac{1}{2}-\epsilon})^{N}$, for any $\epsilon>0$, where the constant $c(\epsilon,f,\alpha)$ depends, aside from $\epsilon$, only on the set of places of $K$ ramifying in $K_{\infty}(f,\alpha)$. This constant can be made explicit.

We next turn to the case of \emph{post-critically infinite} unicritical polynomials, i.e. $f:=x^d-c$ for some $d \in \mathbb{Z}_{\geq 2}$ and $c$ in $K$. In what follows we focus on the number field case, since for function fields one has, in the post-critically infinite case, tools such as Vojta's conjecture. 
\begin{theorem} \label{exponential: unconditional}
Let $K$ be a number field, let $d$ be in $\mathbb{Z}_{\geq 2}$ and $c \in K$. Suppose that $f:=x^d-c$ is post-critically infinite. Then for each $\alpha$ in $K$ there exists a positive constant $c(f,\alpha)$ such that
$$\#G_N(f,\alpha) \geq \textup{exp}(c(f,\alpha) \cdot N), 
$$
for all $N$ in $\mathbb{Z}_{\geq 1}$. 
\end{theorem}
The lower bound in Theorem \ref{exponential: unconditional} is unconditional, however, as explained at the beginning of this introduction, it is far off from the conjectural growth of $D_N(f,\alpha)$ for post-critically infinite $f$. Nevertheless, to the best of our knowledge, this is the first general exponential lower bound for $D_N(f,\alpha)$ proven unconditionally for all unicritical post-critically infinite polynomials over number fields. With some extra work, one can provide the constant $c(f,\alpha)$ explicitly, see the discussion immediately before Proposition \ref{there must be a useful place of periodic reduction}. We remark that the proof of Theorem \ref{exponential: unconditional} adapts very easily to cover also all $x^d+c$ that are PCF and with periodical critical orbit.

As the reader can learn at the end of the proof of Theorem \ref{exponential: unconditional}, the lower bound we prove there is of the form $C \cdot \text{deg}(f)^{\frac{N}{n_0}}$, where $C$ is a positive constant and $n_0$ is the smallest integer such that there exists $v$ a non-archimedean place of $K$ satisfying $v(f^{n_0}(0))>0$ and $v(\alpha+c) \geq 0$. In particular if we take $\alpha:=0$ and $c$ to be not the inverse of an element of $\mathcal{O}_K[\frac{1}{\text{deg}(f)}]$ we get a lower bound of the form $C \cdot \text{deg}(f)^N$. As we next explain, this special case is also a previously known instance of Theorem \ref{exponential: unconditional}, which brings us to the last part of this section, namely the relation between our lower bounds and the notion of stability.  

In case one shows that the pair $(f,\alpha)$ is eventually stable then it follows an exponential lower bound of the form $\#G_N(f,\alpha) \geq C \cdot \text{deg}(f)^N$, for all $N$ in $\mathbb{Z}_{\geq 1}$, where $C$ is a positive constant \cite[Proposition 2.1]{Jones--Levin}. Eventual stability had been established in case $f=x^d+c, \alpha:=0$ and $v(c)>0$ for some non-archimedean place $v$ of $K$ \cite[Theorem 1.6]{Jones--stable} and for $K:=\Q$ and $f:=x^d+c, \alpha:=0$ where $c$ is the inverse of an odd integer different from $-1$ \cite[Corollary 1.5]{Jones--Levin}. From this perspective, Theorem \ref{exponential: unconditional} can be viewed as a generalization \footnote{Within the specialization from eventually stable to exponential lower bounds.} of \cite[Theorem 1.6]{Jones--stable}, where we get an exponential lower bound whose quality is inversely proportional to the amount of time the critical orbit is not divided by some prime where $\alpha$ is integral. Indeed, the proof of Theorem \ref{exponential: unconditional} works more in general for all $f=x^d+c$ where such a valuation $v$ exists and not only in the post-critically infinite case, where such a $v$ can be always found. For more details on how we exploit also orbits of length larger than $1$ (the case of orbits of length $1$ being in particular covered by \cite[Theorem 1.6]{Jones--stable}) we refer to the next subsection, where we give a quick summary of the key steps of the proofs. 

\subsection{Outline of the proofs:}
In the proof of Theorem \ref{exponential:GRH}, the main feature of PCF polynomials that we exploit is that the entire extension $K_{\infty}(f,\alpha)/K$ ramifies only at finitely many places, and, even more, each of the polynomials $f^{N}-\alpha$ has discriminant supported at a fixed finite set of places $S(f,\alpha)$ independent of $N$ (Proposition \ref{only finitely many primes}). This forces the smallest splitting prime of the extension $K_N(f,\alpha)/K$, outside of $S(f,\alpha)$, to have norm at least $\text{deg}(f)^N$ (Proposition \ref{splitting primes are large}). Hence, using GRH (Proposition \ref{GRH: small splitting prime}), one deduces that the logarithm of the absolute discriminant $\text{log}(d_{K_N(f,\alpha)})$ must be exponentially large in $N$ (Proposition \ref{GRH: discriminant is large}). Since $K_N(f,\alpha)$ is ramified only at the fixed finite set $S(f,\alpha)$, this forces $[K_N(f,\alpha):K]$ to be exponentially large in $N$ (Proposition \ref{large discriminant+bounded ramif. locus--->large degree}). Putting everything together then yields Theorem \ref{exponential:GRH}, as we articulate at the end of Section \ref{proof}. We remark that GRH has played a role in previous results in arithmetic dynamics, such as for instance \cite{GRH2} and \cite{GRH1}, in studying the density of prime divisors of linear recurrence sequences. However this assumption here is used in a drastically different manner and for entirely different reasons: to give a sense of the difference we remark that in Theorem \ref{exponential:GRH} the relevant recurrence sequences are supported in finite sets (and the polynomials have degree at least $2$). We also remark that it would be enough to have a real number $\sigma<1$ such that $\text{Re}(s) \leq \sigma $ holds for any $s$ zero of any zeta function. \\
It would be interesting to adapt Theorem \ref{exponential:GRH} to function fields, and obtain unconditional bounds therein. For that one needs to take care of two issues. One comes from constant field extensions, which has serious consequences on the application of Chebotarev. If the constant field extension in $K_N(f,\alpha)$ grows even only linearly in $N$, then the strategy of proof used here does not work. A possible way out could be to establish a dichotomy, which shows that either the constant field grows very slowly in $N$ or exponentially fast, for all PCF polynomials $f$. However, at the moment, we do not know whether this holds and how one could prove it. The other issue comes from wild ramification. Contrarily to a local number field, a local function field can have extensions of arbitrarily large (wild) ramification and given degree, and this will create problems in establishing a good analogue of Proposition \ref{large discriminant+bounded ramif. locus--->large degree}. A way out for the second problem can be found by imposing $\text{deg}(f)<p$, which imposes all the ramification to be tame. This assumption seems particularly natural in the context of PCF polynomials over function fields, see for instance \cite[page 2350]{attracting cycles}.\\
In the proof of Theorem \ref{exponential: unconditional} we explore a generalization of \cite[Lemma 4.1]{B} to general unicritical PCF with \emph{periodic} critical orbit (Proposition \ref{periodical--->root of unity}): \cite{B} focused on the polynomial $x^2-1$, where the critical orbit has period $2$, and showed, among other things, in that case how to explicitly construct larger and larger $2$-power roots of unity by moving $2$ steps at the time up on $T_{\infty}(f,\alpha)$. In Proposition \ref{periodical--->root of unity} we show how to construct larger and larger $d$-power roots of unity by moving $n_0$-steps at the time up on $T_{\infty}(f,\alpha)$, where $n_0$ is the period of the critical orbit. After having done this, we simply observe that a post-critically infinite unicritical polynomial must be PCF with periodic critical orbit modulo some prime that is sufficiently non-degenerate from the standpoint of Proposition \ref{periodical--->root of unity} (this is done in Proposition \ref{there must be a useful place of periodic reduction}). Hence we lower bound $D_N(f,\alpha)$ by the residue field degree at this suitably chosen prime. The conclusion then follows easily. This proof is in line with \cite[Remark 1.4]{local arboreal}, where one uses an auxiliary prime to lower bound the size of global arboreal Galois groups by the size of local arboreal Galois groups. It would be very interesting to devise a way to use several auxiliary periodic primes at the same time, keeping in mind tools such as \cite{Shparlinski} giving some control on the number of prime divisors dividing the critical orbit. To this end \cite{Pink} might provide useful guidance.  
\section*{Acknowledgements}
This project started as a correspondence with Jorge Mello, Alina Ostafe and Igor Shparlinski, who drew my attention to this problem. I warmly thank them for generously sharing with me their own ideas on this subject and for their feedback on previous versions of this work.  
I am grateful to Andrea Ferraguti for his comments on an earlier draft that improved the presentation, for noticing a strenghtening of Theorem \ref{exponential:GRH} and for insightful conversations. I thank Rob Benedetto, Andrew Bridy, Rafe Jones, Valentijn Karemaker and Tom Tucker for valuable feedback. I am grateful to Daniel El-Baz, Pieter Moree and Marco Pagano for giving many suggestions that drastically improved the presentation.  I wish to thank the Max Planck Institute for Mathematics in Bonn for its great working conditions and an inspiring atmosphere.  I gratefully acknowledge the financial support through EPSRC Fellowship EP/P019188/1, ``Cohen--Lenstra heuristics, Brauer relations, and low-dimensional manifolds''.  

\section{Proof of Theorem \ref{exponential:GRH}} \label{proof}
Let $K$ be a number field and $d$ be in $\mathbb{Z}_{\geq 2}$. Until the final paragraph of this section $f$ denotes a PCF polynomial over $K$ of degree $d$ and $\alpha$ denotes any element of $K$ outside of the critical orbits of $f$. At the very end of the section we will reduce the proof of Theorem \ref{exponential:GRH} to this case. Let us begin with a standard fact. 
\begin{proposition} \label{size of preimages}
Let $K,f,\alpha$ be as above. Then
$$\#f^{-N}(\alpha)=d^N,
$$
for all $N$ in $\mathbb{Z}_{\geq 0}$. 
\begin{proof}
This follows immediately, for instance, by induction from \cite[Lemma 2.6]{Jones: LMS}, with $g(x):=x-\alpha$. We now provide a direct argument as well. The desired conclusion is obviously true for $N=0$. Suppose by contradiction that the conclusion does not hold for some positive value of $N$. Let $M$ be the smallest positive integer with $\#f^{-M}(\alpha)<d^M$. Then there must necessarily be $\beta \in f^{-M+1}(\alpha)$ with $f^{-1}(\beta)$ possessing a critical point, $x_0$, of $f$: otherwise each of the $d^{M-1}$ points of $f^{-M+1}(\alpha)$ would have pairwise disjoint $d$-sets of preimages yielding $\#f^{-M}(\alpha)=d^M$, which cannot hold by definition of $M$. Hence $\alpha$ is in the orbit of $x_0$, which means that $\alpha$ is in a critical orbit, which has been excluded from the beginning of this section.  
\end{proof}
\end{proposition}
We now recall how the assumption that $f$ is PCF, along with the assumption that $\alpha$ is outside the critical orbits, forces the various discriminants $\text{Disc}(f^N-\alpha)$ to be supported in a uniform finite set of primes. 
\begin{proposition} \label{only finitely many primes}
 Let $K,f,\alpha$ be as above. Then there exists a finite set $S(f,\alpha)$ of places of $K$, such that $f$ is integral outside of $S(f,\alpha)$ and for any finite prime $\mathfrak{p}$ of $K$ outside of $S(f,\alpha)$ and all positive integers $N$, we have that $v_{\mathfrak{p}}(\textup{Disc}(f^{N}-\alpha))=0$. 
\begin{proof}
This is an immediate consequence of the formula for $\text{Disc}(f^{N}-\alpha)$ given in \cite[Lemma 2.6]{Jones: LMS}, plugging in $g(x):=x-\alpha$ and keeping in mind that Proposition \ref{size of preimages} forces each factor in that formula to be non-zero. 
\end{proof} 
\end{proposition} 
For the rest of this section we fix such a finite set of places $S(f,\alpha)$ provided by Proposition \ref{only finitely many primes}. 

For a number field $L/\Q$, we denote by $d_L$ the absolute discriminant of $L$. Fix now $\epsilon \in \mathbb{R}_{>0}$ any positive constant. The following result is a consequence of the work of Lagarias--Odlyzko \cite{OL}. 
\begin{proposition} \label{GRH: small splitting prime}
Assume GRH. Then there exists a positive constant $c_1(f,\alpha)$ such that for each $N$ in $\mathbb{Z}_{\geq 1}$ there exists a prime $\mathfrak{p}_N$ in $K$ outside of $S(f,\alpha)$ with the following two properties. \\
$(1)$ $\mathfrak{p}_N$ splits completely in $K_N(f,\alpha)$. \\
$(2)$ $\#\mathcal{O}_K/\mathfrak{p}_N \leq c_1(f,\alpha,\epsilon) \cdot \textup{log}(d_{K_N(f,\alpha)})^{2+\epsilon}$
\begin{proof}
We apply \cite[Theorem 4]{Serre} with the same $K$, with $E:=K_N(f,\alpha)$, with 
$$x:=c_1(f,\alpha,\epsilon) \cdot \text{log}(d_E)^{2+\epsilon},
$$
and with $C$ the conjugacy class of the identity. The constant $c_1(f,\alpha)$ will be chosen sufficiently large, as we explain below. Then \cite[Theorem 4]{Serre} combined with the negation of the proposition yields that $\text{Li}(c_1(f,\alpha, \epsilon)\text{log}(d_E)^{2+\epsilon})$ is at most
$$c_6 \text{log}(d_E)^{1+\frac{\epsilon}{2}}\cdot (\text{log}(d_E)+[E:\Q]\text{log}(c_1(f,\alpha,\epsilon)\text{log}(d_E)^{2+\epsilon}))+\#S(f,\alpha). 
$$  
Here $c_6$ is the positive absolute constant taken from the statement of \cite[Theorem 4]{Serre}. 

Plugging in the Minkowski's bound for $[E:\Q]$, we conclude in particular that the quantity $\text{Li}(c_1(f,\alpha,\epsilon)\text{log}(d_E)^{2+\epsilon})$ is at most 
$$ \text{log}(d_E)^{1+\frac{\epsilon}{2}}\cdot (\text{log}(d_E)+a_{\text{abs}} \cdot \text{log}(d_E) \cdot \text{log}(c_1(f,\alpha,\epsilon)\text{log}(d_E)^{2+\epsilon}))+\#S(f,\alpha),
$$ 
for some absolute constant $a_{\text{abs}}$. Now we see that it is possible to chose $c_1(f,\alpha,\epsilon)$ so large that the last inequality cannot hold for any $E$. To avoid this contradiction, that arose from the negation of the present statement, we obtain the desired conclusion. 
\end{proof}
\end{proposition}
On the other hand we have the following elementary fact.
\begin{proposition} \label{splitting primes are large}
Let $\mathfrak{p}$ be a finite place of $K$ outside of $S(f,\alpha)$ and let $N$ be in $\mathbb{Z}_{\geq 0}$. Suppose that $\mathfrak{p}$ splits completely in $K_N(f,\alpha)/K$. Then
$$\#(\mathcal{O}_K/\mathfrak{p}) \geq d^N. 
$$
\begin{proof}
Take $\mathfrak{p}$ as in the statement. Let $\mathbb{F}_\mathfrak{p}:=\mathcal{O}_K/\mathfrak{p}$ be the residue field of $\mathfrak{p}$. Choose any prime $\widetilde{\mathfrak{p}}$ above $\mathfrak{p}$ in $K_N(f,\alpha)$. Since $\mathfrak{p}$ splits completely, we have that the natural inclusion map induces an equality 
$$\mathcal{O}_{K_N(f,\alpha)}/\widetilde{\mathfrak{p}}=\mathbb{F}_{\mathfrak{p}}.
$$
Observe that, by definition of $S(f,\alpha)$, the set $f^{-N}(\alpha)$ consists of integral elements locally at $\widetilde{\mathfrak{p}}$ and thus we can reduce each of its elements modulo $\widetilde{\mathfrak{p}}$. 

When reduced modulo $\widetilde{\mathfrak{p}}$, the various elements of $f^{-N}(\alpha)$ need to be distinct in view of Proposition \ref{only finitely many primes}. It follows that
$$\#\mathbb{F}_{\mathfrak{p}}=\#\mathcal{O}_{K_N(f,\alpha)}/\widetilde{\mathfrak{p}} \geq \#f^{-N}(\alpha)=d^N,
$$
where the last equality follows from Proposition \ref{size of preimages}. This is precisely the desired inequality. 
\end{proof}
\end{proposition}
Combining Proposition \ref{GRH: small splitting prime} and Proposition \ref{splitting primes are large}, we conclude the following. 
\begin{corollary} \label{GRH: discriminant is large}
Assume GRH. Then there exists a positive constant constant $c_2(f,\alpha,\epsilon)$ such that for all $N$ in $\mathbb{Z}_{\geq 1}$, we have that
$$\textup{log}(d_{K_N(f,\alpha)}) \geq c_2(f,\alpha,\epsilon) \cdot (d^{\frac{1}{2+\epsilon}})^N.  
$$
\end{corollary}
Our next step lower bounds the degree in terms of the discriminant for extensions ramifying only above a fixed set of places.   
\begin{proposition} \label{large discriminant+bounded ramif. locus--->large degree}
Let $F$ be a number field, and let $S$ be a finite set of finite places in $F$. Then there is a positive constant $c(S)$ such that for any extension $E/F$ unramified at all finite places outside of $S$, we have that
$$c(S) \cdot [E:\Q]\textup{log}([E:\Q]) \geq \textup{log}(d_E). 
$$
\begin{proof}
We will show that this follows from \cite[Proposition 4']{Serre}. We first explain the correspondence between the notation therein and the notation of this proposition. The $K$ in that proposition is our $F$, the notation for $E$ agrees. The set $P(E/K)$ in \cite{Serre} consists for us of the set of rational primes below a prime of $F$ that is in $S$. The number $n$ is $[E:F]$ and the number $n_E$ is $[E:\Q]$. Hence the quantity $(1-\frac{1}{n}) \cdot \sum_{p \in P(E/K)}\text{log}(p)$ is a nonnegative constant, which we denote by $c_1(S)$, and likewise $\#P(E/K)$ and $\text{log}(d_K)$ are other nonnegative constants, which we denote by $c_2(S)$ and $c_3(S)$, respectively. It now follows from \cite[Proposition 4']{Serre} that
$$\frac{\text{log}(d_E)}{[E:\Q]} \leq  c_3(S)+c_1(S)+c_2(S)\text{log}([E:\Q]),
$$
which implies the desired conclusion. 

\end{proof}
\end{proposition}

We are now ready to prove Theorem \ref{exponential:GRH}. \\
\emph{Proof of Theorem \ref{exponential:GRH}:} We keep assuming that $\alpha$ is outside of all the critical orbits of $f$, as we did so far in this section, and we prove the theorem under this assumption. We explain only at the end how to reduce to this case when we only assume that $(f,\alpha)$ is not exceptional. Thanks to Proposition \ref{large discriminant+bounded ramif. locus--->large degree} and Proposition \ref{only finitely many primes} we know that for each $\epsilon>0$ we can write
$$[K_N(f,\alpha):K] \geq c'(f,\alpha, \epsilon)\text{log}(d_{K_N(f,\alpha)})^{\frac{1}{1+\epsilon}},
$$
for some positive constant $c'(f,\alpha, \epsilon)$. By Corollary \ref{GRH: discriminant is large} we conclude that there exists a positive constant $\widetilde{c}(\epsilon,f,\alpha)$ such that
$$[K_N(f,\alpha):K] \geq \widetilde{c}(\epsilon,f,\alpha) \cdot (d^{\frac{1}{1+\epsilon} \cdot \frac{1}{2+\epsilon}})^N,
$$
for all $N$ in $\mathbb{Z}_{\geq 1}$. Since $\epsilon>0$ is arbitrary, we can rewrite this by saying that for each $\epsilon>0$ there exists a positive constant $c(f,\alpha, \epsilon)$ such that
$$[K_N(f,\alpha):K] \geq c(f,\alpha, \epsilon) \cdot (d^{\frac{1}{2}-\epsilon})^N,
$$
for all $N$ in $\mathbb{Z}_{\geq 1}$. This gives in particular the desired conclusion under the assumption that $\alpha$ is outside of the critical orbit of $f$. We now drop this assumption and we only assume that $(f,\alpha)$ is not exceptional. This will conclude the proof in complete generality. Since $(f,\alpha)$ is not exceptional we have that $T_{\infty}(f,\alpha)$ is infinite, and thus, since $f$ is PCF, there must be $\beta$ in $T_{\infty}(f,\alpha)$ such that $\beta$ is outside of all the critical orbits of $f$. Then the above proof gives the conclusion for $D_N(f,\beta)$, where the number field $K$ is replaced with the number field $K(\beta)$. This clearly implies the desired lower bound for $D_N(f,\alpha)$.   

\section{Proof of Theorem \ref{exponential: unconditional}}
Our main tool will be a construction of roots of unity in the arboreal fields of a PCF unicritical polynomial $f:=x^d-c$ over a field $F$ of characteristic coprime with $d$, under the assumption that the critical orbit, which equals to the orbit of $0$, is periodic under the action of $f$. This construction is inspired by \cite[Lemma 4.1]{B}, which dealt with the special case $f:=x^2-1$, where the length of the period is $2$. We have the following.  
\begin{proposition} \label{periodical--->root of unity}
Let $d$ be an integer, let $F$ be a field of characteristic coprime to $d$. Let $c$ be in $F$. Suppose there is a positive integer $n_0$ such that $f^{n_0}(0)=0$, where $f:=x^d-c$. Let $\alpha$ be an element of $K$ different from $-c$.

Then, for each positive integer $N$, the field $F_N(f,\alpha)$ contains an element of multiplicative order equal to 
$$d^{\textup{max}(\lfloor \frac{N}{n_0} \rfloor -1,0)+1}
$$ 
\begin{proof}
Since $\alpha$ is different from $-c$ and since $d$ is coprime to the characteristic of $F$, we know that there are $d$ distinct non-zero elements in $f^{-1}(\alpha)$, that form a free transitive set under the action of the cyclic group of order $d$ given by
$$\mu_d(F_1(f,\alpha)):=\{\zeta \in F_1(f,\alpha): \zeta^d=1\}. 
$$ 
Take $\zeta_d$ a generator of $\mu_d(F_1(f,\alpha))$, pick an element $\beta$ in $f^{-1}(\alpha)$, and denote by $\beta'$ the element $\zeta_d \cdot \beta$. 

Observe that the polynomial $f^{n_0}-\beta$  is a polynomial in $x^d$ of degree $d^{n_0}$. It follows that we can pick a vector $(\gamma_1(1), \ldots,\gamma_{d^{n_0-1}}(1))$ in $(f^{-n_0}(\beta))^{d^{n_0-1}}$ such that
$$f^{n_0}-\beta=(x^d-\gamma_1(1)^d) \ldots (x^d-\gamma_{d^{n_0-1}}(1)^d).
$$
We now evaluate in $0$ and take advantage of the fact that $f^{n_0}(0)=0$ to obtain that
$$-\beta=(-1)^{d^{n_0-1}}\cdot (\gamma_1(1) \ldots \gamma_{d^{n_0}-1}(1))^d.
$$
We can do the same for $\beta'$ and obtain that
$$-\beta'=(-1)^{d^{n_0-1}}\cdot (\gamma_1(1)' \ldots \gamma_{d^{n_0}-1}(1)')^d,
$$
with $(\gamma_1(1)', \ldots,\gamma_{d^{n_0-1}}(1)')$ in $(f^{-n_0}(\beta'))^{d^{n_0-1}}$.

Now we can iterate the procedure above and find a function $g:\mathbb{Z}_{\geq 1} \to \mathbb{F}_2$ such that for each $i$ in $\mathbb{Z}_{\geq 1}$ there are two vectors $(\gamma_1(i), \ldots, \gamma_{d^{i \cdot n_0-i}}(i)) \in (f^{-n_0 \cdot i}(\beta))^{d^{i \cdot n_0-i}}$ with
$$-\beta=(-1)^{g(i)} \cdot (\gamma_1(i) \ldots \gamma_{d^{i \cdot n_0-i}}(i))^{d^{i}}
$$
and
$$-\beta'=(-1)^{g(i)} \cdot (\gamma_1(i)' \ldots \gamma_{d^{i \cdot n_0-i}}(i)')^{d^{i}}
$$
Recalling that the elements $\beta$ and $\beta'$ are non-zero and their ratio equals $\zeta_d$, we obtain that, for each $i$ in $\mathbb{Z}_{\geq 1}$, $\zeta_d$ is a $d^{i}$-th power in $K_{i \cdot n_0+1}(f,\alpha)$. That means that for each $i$ in $\mathbb{Z}_{\geq 0}$, one has inside the field $K_{i \cdot n_0+1}(f,\alpha)$ an element of multiplicative order precisely equal to $d^{i+1}$. This is the desired conclusion. 
%Now observe that none of the elements in $f^{-n_0}(\beta),f^{-n_0}(\beta')$ can possibly be $0$. Indeed, since $f^{n_0}(0)=0$, this would force either $\beta$ or $\beta'$ to be $0$, which would imply that $\alpha=-c$, which has been assumed not to be the case.  

%Now on the multi-sets of roots of $f^{n_0}(x)-\beta$ and $f^{n_0}(x)-\beta'$, respectively, counted with multiplicity, there is a natural action of $\mu_d(F_1(f,\alpha))$, which is free and transitive. Call $\mathcal{S}, \mathcal{S}'$ a set of representatives for each orbit. Then 
%$$(-1)^{\#\mathcal{S}}(\prod_{\gamma \in \mathcal{S}}\gamma)^d$$
%equals
%$$\prod_{\tau \in f^{-n_0+1}(\beta)}(-c-\tau)=f^{n_0-1}(-c)-\beta=f^{n_0}(0)-\beta=-\beta, 
%$$
%and 
%$$(-1)^{\#\mathcal{S}'}(\prod_{\gamma' \in \mathcal{S}'}\gamma')^d$$
%equals
%$$\prod_{\tau' \in f^{-n_0+1}(\beta')}(-c-\tau')=f^{n_0-1}(-c)-\beta'=f^{n_0}(0)-\beta=-\beta'.
%$$
%Recalling that $\beta,\beta'$ are non-zero and taking the ratio of these two expressions, we conclude that $\zeta_d$ is a $d$-th power in $F_{n_0+1}(f,\alpha)$. Now, since $f^{-n_0}(\beta),f^{-n_0}(\beta')$ consists entirely of non-zero elements, as we argued above, we can repeat the same procedure on each $\gamma$ and find out that $\zeta_d$ is $d^2$ power in $F_{2n_0+1}$. Iterating we find that $\zeta_d$ is a $d^k$-power in $F_{kn_0+1}(f,\alpha)$, which means that $F_{kn_0+1}(f,\alpha)$ has an element of multiplicative order precisely equal to $k+1$. This is precisely the desired conclusion. 
\end{proof}
\end{proposition}
The next proposition will offer us a useful prime of periodic reduction in the post-critically infinite case as an immediate consequence of a result in \cite{Tucker}, which, as explained therein, is a consequence of Siegel's theorem on integral points. As an alternative, we remark that in this special case one could easily adapt the arguments given in \cite[Thm 2.6]{FP}. We remark that \cite[T.6]{Shparlinski} obtains a much stronger conclusion for polynomials over $\Q$. 
\begin{proposition}\label{there must be a useful place of periodic reduction}
Let $K$ be a number field, $d$ in $\mathbb{Z}_{\geq 2}$ and $c \in K$ such that $x^d-c$ is post-critically infinite. Let $\alpha$ be in $K$ different from $-c$. Then there exists a positive integer $n_0$ and a prime $\mathfrak{p}$, coprime to $d$, with $v_{\mathfrak{p}}(\alpha+c)=0$ and such that
$$v_{\mathfrak{p}}(f^{n_0}(0))>0. 
$$
\begin{proof}
We will show that this is an immediate consequence of \cite[Proposition 1.6, (a)]{Tucker}, which we can apply since $c \neq 0$. The orbit of $0$ takes negative valuation only at those places where $c$ does, so finitely many. Since $\alpha \neq c$ there are only finitely many valuation where $\alpha+c$ is non-zero, and the places above prime divisors of $d$ are finite. Call $S$ the finite set obtained collecting all these places listed so far. By \cite[Proposition 1.6, (a)]{Tucker}, we know that there exists $n_0$ in $\mathbb{Z}_{\geq 1}$ such that $f^{n_0}(0)$ is not an $S$-unit. This means that there exists a place $\mathfrak{p}$ outside of $S$ such that $v_{\mathfrak{p}}(f^{n_0}(0)) \neq 0$. Since $\mathfrak{p}$ is not in $S$, the only possibility is that $v_{\mathfrak{p}}(f^{n_0}(0))>0$. Furthermore since $\mathfrak{p}$ is not in $S$, we have that $\mathfrak{p}$ is coprime to $d$ and $v_{\mathfrak{p}}(\alpha +c)=0$. This is precisely the desired conclusion. 
\end{proof}
\end{proposition}
To turn the combination of Proposition \ref{periodical--->root of unity} and Proposition \ref{there must be a useful place of periodic reduction} into a lower bound we need the following basic fact. 
\begin{proposition} \label{cyclotomic extension of finite fields}
Let $\mathbb{F}$ be a finite field of cardinality coprime to an integer $d$. Then there is a positive constant $c(d,\mathbb{F})$ such that the following holds. If $\mathbb{L}/\mathbb{F}$ is a finite extension such that $\mathbb{L}$ has an element of multiplicative order equal to $d^s$, then
$$[\mathbb{L}:\mathbb{F}] \geq  c(d, \mathbb{F}) \cdot d^s.
$$
\begin{proof}
Denote by $q:=\#\mathbb{F}$. Denote by $\Gamma$ the closure of the group generated by $q$ in $\prod_{l|d}\mathbb{Z}_l^{*}$. Observe that $\Gamma$ is open in $\prod_{l|d}\mathbb{Z}_l^{*}$. Define
$$U(s):=\{x \in \prod_{l|d}\mathbb{Z}_l^{*}: x \equiv 1 \ \text{mod} \ d^s\}. 
$$
In order to prove the statement we need to show that there exists a positive constant $c(d,\mathbb{F})$ such that for each positive integer $s$ we have that
$$[\Gamma:\Gamma \cap U(s)] \geq c(d,\mathbb{F}) \cdot d^s.
$$
Indeed for $\mathbb{L}$ as in the statement we have
$$[\mathbb{L}:\mathbb{F}] \geq [\Gamma:\Gamma \cap U(s)] .
$$
Observe that there exists an isomorphism of profinite group, $\phi_d$, sending
$$\phi_d: U(2) \to \prod_{l|d}\mathbb{Z}_l,
$$
such that $\phi_d(U(s))=d^{s-2} \cdot \prod_{l|d}\mathbb{Z}_l$ for each $s$ at least $2$. Finally observe that 
$$\phi_d(\Gamma \cap U(2))=\prod_{l|d}l^{f_l(\mathbb{F})}\mathbb{Z}_l,
$$ 
for suitable positive integers $f_l(\mathbb{F})$. This is an immediate consequence of the fact that $\Gamma$ is open in $\prod_{l|d}\mathbb{Z}_l^{*}$. Define $m_0(\mathbb{F})$ to be the maximum of the various $f_l(\mathbb{F})$, with $l|d$. We now see that
$$[\Gamma \cap U(2): \Gamma \cap U(s)] \geq d^{s-m_0(\mathbb{F})},
$$
for each $s \geq m_0(\mathbb{F})$, which gives the desired conclusion. 

\end{proof}
\end{proposition}
We are now ready to prove Theorem \ref{exponential: unconditional}. \\
\emph{Proof of Theorem \ref{exponential: unconditional}:} Assume first that $\alpha$ is different from $-c$. Thanks to Proposition \ref{there must be a useful place of periodic reduction}, we can find a prime $\mathfrak{p}$ such that $f$ reduces modulo $\mathfrak{p}$ to a polynomial of the form $x^d-c$ with the orbit of $0$ periodic of period $n_0>0$, with $v_{\mathfrak{p}}(\alpha+c)=0$ and with $d$ coprime to $\mathfrak{p}$. Let us fix such a place $\mathfrak{p}$. Let now $N$ be a positive integer and take $\widetilde{\mathfrak{p}}_N$ a place above $\mathfrak{p}$ in $K_N(f,\alpha)$. Let us denote by $\mathbb{F}_{\widetilde{\mathfrak{p}}_N}$ and $\mathbb{F}_{\mathfrak{p}}$ the corresponding residue fields. Observe that $c$ is integral locally at $\mathfrak{p}$, since $v_{\mathfrak{p}}(f^{n_0}(0))>0$, while $v_{\mathfrak{p}}(c)<0$ certainly implies $v_{\mathfrak{p}}(f^{n_0}(0))<0$. Likewise, we know that $\alpha$ is integral locally at $\mathfrak{p}$, since, if not, we could not have $v_{\mathfrak{p}}(\alpha+c)=0$, since $v_{\mathfrak{p}}(c) \geq 0$. Therefore we can reduce the entire set $f^{-N}(\alpha)$ modulo $\widetilde{\mathfrak{p}_N}$ and apply Proposition \ref{periodical--->root of unity}. Hence we deduce that there is in $\mathbb{F}_{\widetilde{\mathfrak{p}}_N}$ an element of multiplicative order larger than $d^{\frac{N}{n_0}-1}$. Hence, via Proposition \ref{cyclotomic extension of finite fields}, we deduce that
$$[\mathbb{F}_{\widetilde{\mathfrak{p}}_N}:\mathbb{F}_{\mathfrak{p}}] \geq c(d,\mathfrak{p}) \cdot d^{\frac{N}{n_0}},  
$$  
where $c(d,\mathfrak{p})$ is a positive constant depending only on $d$ and the place $\mathfrak{p}$. Recalling that $[\mathbb{F}_{\widetilde{\mathfrak{p}}_N}:\mathbb{F}_{\mathfrak{p}}]$ divides $[K_N(f,\alpha):K],
$
we obtain the desired conclusion. We are now left with the case that $\alpha=-c$, which is equivalent to $f^{-1}(\alpha)=\{0\}$. We can redo exactly the same argument for $\alpha':=0$, which is certainly different from $-c$ since $f$ is post-critically infinite. This provides a satisfactory bound for $[K_N(f,\alpha'):K]=[K_{N+1}(f,\alpha):K]$ and thus for $D_N(f,\alpha)$ as desired.


\begin{thebibliography}{99}

\bibitem{B}
F. Ahmad, R. Benedetto, J. Cain, G. Carroll, L. Fang, ``The arithmetic Basilica: a quadratic PCF arboreal Galois group", \textit{arXiv preprint: 1909.0039}.

\bibitem{local arboreal}
J. Anderson, S. Hamblen, B. Poonen, L. Walton, ``Local arboreal representations", IMRN, Volume 2018, Issue 19, October 2018, Pages 5974–5994.

\bibitem{Andrews--Petsche}
J. Andrews, C. Petsche, ``Abelian extensions in dynamical Galois theory" Algebra and
Number Theory, 14(7): 1981-1999 (2020).

\bibitem{Benedetto: surjectivity}
R. Benedetto, J. Juul, ``Odoni's conjecture for number fields" Bulletin of the London Mathematical Society 51, $\#2$ (2019), 237-250.

\bibitem{attracting cycles}
R. Benedetto, P. Ingram, R. Jones, A. Levy, ``Attracting cycles in $p$-adic dynamics and height bounds for post-critically finite maps", Duke Math. J. 163 (13) 2325 - 2356, 1 October 2014.

\bibitem{BIJMST}
R. Benedetto, P. Ingram, R. Jones, M. Manes, J.H. Silverman, and T.J. Tucker, 
``Current trends and open problem in arithmetic dynamics",  Bull. Amer. Math. Soc., {\bf 56} (2019), 611--685.
%
%\bibitem{BIJJLMRS}
%A. Bridy, P. Ingram, R. Jones, J. Juul, A. Levy, M. Manes, S. Rubinstein-Salzedo, and J. H. Silverman,  
%`Finite ramification for preimage fields of post-critically finite morphisms',
%  {\it Math. Res. Lett.}, {\bf 24} (2017), no. 6, 1633-1647.

\bibitem{Bush--Hindes--Looper}
M. Bush, W. Hindes, N. Looper, ``Galois groups of iterates of some unicritical polynomials" Acta Arithmetica 181 (2017): 57-73. 

\bibitem{Cremona}
J. E. Cremona, ``On the Galois groups of the iterates of $x^2+1$", Mathematika, 36(2):259–261 (1990), 1989. 
%\bibitem{Dub}
%A. Dubickas, \textit{Algebraic numbers with bounded degree and Weil height}, 
%Bull. Aust. Math. Soc. 98 (2018), 212--220.
\bibitem{FM}
A. Ferraguti, G. Micheli, ``An equivariant isomorphism theorem for mod $\mathfrak{p}$ reductions of arboreal Galois representations", Transactions of the American Mathematical Society , 373 (2020), no. 12, 8525-8542.
\bibitem{FP}
A. Ferraguti, C. Pagano, ``Constraining images of quadratic arboreal representations", IMRN, Volume 2020, Issue 22, November 2020, Pages 8486-8510.  
\bibitem{FPC}
A. Ferraguti, C. Pagano, D. Casazza, ``The inverse problem for arboreal Galois representations of index two", \textit{arXiv preprint: 1907.08608}.
\bibitem{Gratton--Nguyen--Tucker}
C. Gratton, K. Nguyen, T. Tucker, ``ABC implies primitive prime divisors in arithmetic dynamics", Bulletin of the LMS (2013).
\bibitem{Jones--stable}
S. Hamblen, R. Jones, K. Madhu,  ``The density of primes in orbits of $z^d + c$", Int. Math. Res. Not. 2015(7), 1924-1958.
\bibitem{Hindes}
W. Hindes, ``The Vojta conjecture implies Galois rigidity in dynamical families", Proc. Amer. Math. Soc. 144 (2016), 1931-1938.


          
          
\bibitem{Jones: LMS} R. Jones, `` The density of prime divisors in the arithmetic dynamics of quadratic polynomials", J. London Math. Soc., 2008, 523-544.       

\bibitem{Jones} R. Jones, ``Galois representations from pre-image trees: An arboreal
            survey",  Publ. Math. Besan\c con Alg\`ebre Th\'eorie Nombr.\ , 
            Univ. Franche-Comt{\'e}, Besan{\c c}on, 2013, 107--136. 
\bibitem{Jones--Levin} 
R. Jones, A. Levy, ``Eventually stable rational functions" Int. J. Number Theory 13(9) (2017), 2299-2318.

\bibitem{Win}
J. Juul, H. Krieger, N. Looper, M. Manes, B. Thompson, and L. Walton ``Arboreal representations for rational maps with few critical points", Proceedings of the Women in Numbers 4, \emph{to appear}. 

\bibitem{Kadet}
B. Kadets, ``Large arboreal Galois representations", Journal of Number theory, Volume 210, May 2020, Pages 416-430. 
\bibitem{Tucker}
H. Krieger, A. Levin, Z. Scherr, T. Tucker, Y. Yasufunku, M. Zieve, ``Uniform boundedness of $S$-units in arithmetic dynamics", Pacific J. Math. \textbf{274} (2015), no.2, 97-106. 



\bibitem{GRH2}
J.C. Lagarias, ``The set of primes dividing the Lucas numbers has density 2/3", Pacific J. Math. \textbf{118} (1985) 449–461 (Errata, Pacific J. Math. \textbf{162} (1994), 393-397).
          
\bibitem{OL}
J.C. Lagarias, A.M. Odlyzko,  ``Effective Versions of the Chebotarev Density Theorem" In:
Algebraic Number Fields, L-Functions and Galois Properties (A. Frohlich, ed.), pp. 409-464. New York, London: Academic Press 1977.
\bibitem{Looper}
N. Looper, ``Dynamical Galois groups of trinomials and Odoni's Conjecture" Bulletin of the LMS, \emph{to appear}. 

%
%\bibitem{Mi}  J.S.Milne, `Algebraic number Theory`', Online: \url{www.jmilne.org/math/CourseNotes/}.

%\bibitem{Mur}  R. Murty,  `Prime numbers and Irreducible polynomials', 
%\textit{Amer. Math. Monthly} {\bf 109} (2002), 452--458.
\bibitem{GRH1}
P. Moree, P. Stevenhagen, ``A two-variable Artin conjecture", J. Number Theory 85 (2000) 291–304.

\bibitem{Odoni1}
R. W. K. Odoni. ``The Galois theory of iterates and composites of polynomials", Proc. London
Math. Soc. (3), 51(3):385–414, 1985.
\bibitem{Odoni2}
R. W. K. Odoni. ``On the prime divisors of the sequence $w_n+1 = 1 +w_1 \ldots w_n$", J. London Math.Soc. (2), 32(1):1–11, 1985.



\bibitem{Pink}
R. Pink ``Profinite iterated monodromy groups arising from quadratic polynomials", \textit{arXiv preprint: 1307.5678}.

\bibitem{Serre}
J-P. Serre, ``Quelques applications du th\'eor\`eme de densit\'e de Chebotarev", Publications math\'ematiques de l'I.H.\'E.S., tome 54 (1981), p. 123-201.
\bibitem{Shparlinski}
I. E. Shparlinski, ``Number of different prime divisors of recurrence sequences", Mat. Zametki, 1987, Volume 42, Issue 4, Pages 494-506.  
\bibitem{Silv07}
J. H. Silverman, \textit{The arithmetic of dynamical systems},
Springer-Verlag, New York, 2007.
\bibitem{Spechter}
J. Specter, ``Polynomials with surjective Galois representations exist in every degree", \textit{arXiv preprint: 1803.00434}.

\bibitem{Stoll}
M. Stoll, ``Galois groups over $\Q$ of some iterated polynomials", Arch. Math (Basel), 59(3): 239-244, 1992. 

 

\end{thebibliography}
\end{document}